\documentclass[oneside,10pt]{article}
\usepackage[b5paper]{geometry}	

\usepackage{amsfonts,amsmath,latexsym,amssymb}
\usepackage{theorem}
\usepackage{mathrsfs,upref}
\usepackage{mathptmx}		
	
\usepackage{new}	

\usepackage[dvipdf]{graphicx}

\newtheorem{theorem}{Theorem}
\newtheorem{lemma}{Lemma}

\newtheorem{statement}{Statement}

\theoremstyle{definition}

\newtheorem{remark}{Remark}

\newfont{\eurorm}{eurm10 scaled 1000} \def\eu#1{\mbox{\eurorm #1}}    %
\newfont{\newit}{cmfi10 scaled 1000}  \def\pp#1{\mbox{\newit #1}}     %
\begin{document}

\title[Short title]{Barrow's Inequality and Signed Angle Bisectors}

\author{Branko Male\v sevi\' c, Maja Petrovi\' c}

\address{Branko Male\v sevi\' c, Faculty of Electrical Engineering, University of Belgrade, \\
Bulevar Kralja Aleksandra 73, 11000 Belgrade, Serbia \\
\email{malesevic@etf.rs}}

\address{Maja Petrovi\' c, Faculty of Transport and Traffic Engineering, University of Belgrade, \\
Vojvode Stepe 305, 11000 Belgrade, Serbia \\
\email{majapet@sf.bg.ac.rs}}

\CorrespondingAuthor{Branko Male\v sevi\' c}

\date{19.04.2013}



\thanks{
Research is partially supported by the Ministry of Science and Education of the Republic of Serbia,
Grant No.$\,$III~44006 and ON~174032
}

\begin{abstract}
In this paper we give one extension of {\sc Barrow}'s type inequality in the plane of the triangle
$\triangle\,ABC$ introduce signed angle bisectors.
\end{abstract}

\maketitle

\section{Introduction}

Let triangle $\triangle\,ABC$ be given in Euclidean plane. Denote by $R_A, \, R_B$ and $R_C$
the distances from the arbitrary point $M$ in the plane of $\triangle\,ABC$ to the vertices
$A$, $B$ and $C$ respectively, and denote by  $\ell_a = |MA'|$, $\ell_b  = |MB'|$ and
$\ell_c  = |MC'|$ the length of angle bisectors of $\angle BMC$, $\angle CMA$ and
$\angle AMB$ from the point $M$ respectively (Fig. 1).

\smallskip
{\sc Barrow}'s inequality \cite{EMB}:
\begin{equation}
    \label{E1}
    R_A + R_B + R_C\ge
    2\left(\ell_a + \ell_b + \ell_c\right)
    \smallskip
\end{equation}
is true when $M$ is arbitrary point in the interior of triangle $\triangle\,ABC$. The equality holds
iff triangle $ABC$ is equilateral and point $M$ is its circumcenter. In this paper we consider
a~{\sc Barrow}'s type inequality when $M$ is arbitrary point in the plane of the triangle
$\triangle\,ABC$ introduce signed angle bisectors. Let us notice that inequalities with angle
bisectors recently are considered in papers \cite{Indika}, \cite{Jiang}, \cite{JMI}, \cite{JIPAM}.

\begin{center} 

\vspace*{17.5mm} \hspace*{-25.0mm} \includegraphics*[height=20.0mm,keepaspectratio=true,scale=0.6]{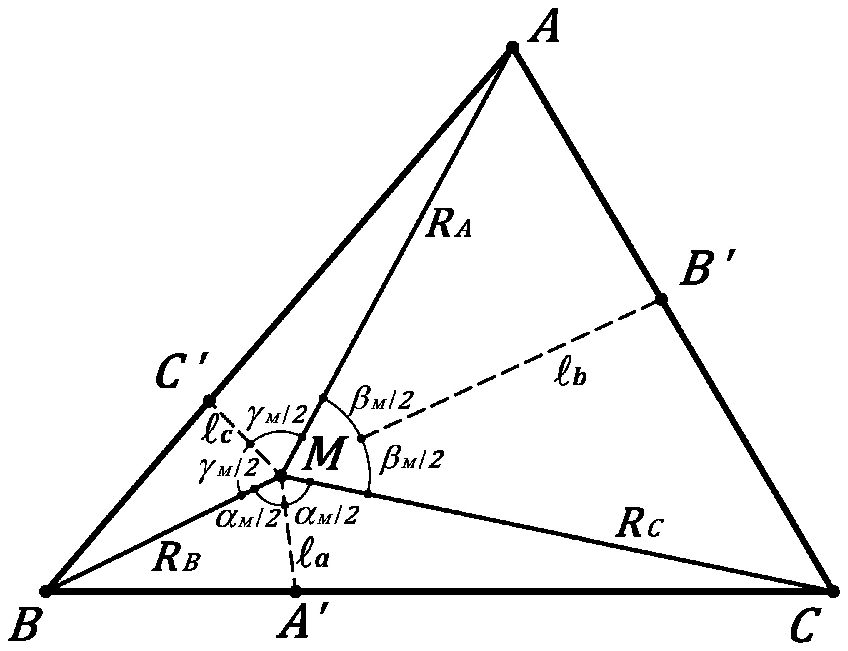}

\smallskip

\noindent
\textit{Figure 1: Barrow's inequality (point $M$ into $\triangle\,ABC$)}
\end{center}

\break

Inequality of {\sc Erd\" os-Mordell} \cite{Erdos}:
\begin{equation}
    \label{E2}
    R_A + R_B + R_C\ge
    2\left(r_a + r_b + r_c\right)
    \smallskip
\end{equation}
is a consequence of inequality of {\sc Barrow}, where $r_a$, $r_b$ and $r_c$ are distances of interior
point $M$ of triangle to the sides $BC$, $CA$ and $AB$ respectively.

\smallskip
Let us notice that topic of the {\sc Erd\" os-Mordell} inequality is current, as it has been shown in recent papers.
{\sc V. Pambuccian} proved that, in the plane of absolute geometry, the {\sc Erd\" os-Mordell} inequality is
an equivalent to non-positive curvature \cite {Pambuccian}. In the paper \cite{MIA} is given an extension of
the {\sc Erd\" os-Mordell} inequality on the interior of the {\sc Erd\"os-Mordell} curve. In relation to
the {\sc Erd\" os-Mordell} inequality {\sc N.~Dergiades} in the paper \cite {Derg} proved one extension of
the {\sc Erd\" os-Mordell} type inequality
\begin{equation}
    \label{E3}
    R_A + R_B + R_C\ge
    \left(\frac{c}{b} + \frac{b}{c}\right)r_a' +
    \left(\frac{c}{a} + \frac{a}{c}\right)r_b' +
    \left(\frac{a}{b} + \frac{b}{a}\right)r_c'
    \smallskip
\end{equation}
where $r_a'$, $r_b'$ and $r_c'$ are signed distances of arbitrary point $M$ in the plane triangle
to the sides $BC$, $CA$ and $AB$ respectively.

\section{The Main Results}

Proof of {\sc Barrow}'s inequality in the paper of {\sc Z. Lu} \cite {Lu} is based on the next theorem.
\begin{statement}
Let $p, \, q, \, r \ge 0$ and $\alpha + \beta + \gamma = \pi $. Then we have the inequality
\begin{equation}
    \label{E4}
    p+q+r\ge
    2\sqrt{qr} \, \cos \alpha \ +
    2\sqrt{pr} \, \cos \beta  \ +
    2\sqrt{pq} \, \cos \gamma.
\end{equation}
\end{statement}

Peculiarity of {\sc Barrow}'s and {\sc Lu}'s proofs are, that is, primarily algebraic. In
{\sc Lu}'s proof, {\sc Barrow}'s inequality follows from positivity of quadratic function
$f\left(x\right) = x^{2} - 2{\big(} \sqrt{r} \, \cos \beta + \sqrt{q}$ $\cos \gamma{\big)}\,
x + q + r - 2\sqrt{qr} \, \cos \alpha$ in the point $x=\sqrt{p}$ with an appropriate geometric
interpretation for $p$, $q$, $r$ and $\alpha$, $\beta$, $\gamma$ (for details see \cite{Lu}).

\smallskip
In this paper we also give one algebraic proof with geometric interpretation for points outside of the triangle
$\triangle\,ABC$. The following theorems are true.

\begin{statement}
\label{St2}
Let $p, \, q, \, r \ge 0$ and $\alpha = \beta + \gamma$. Then we have the inequality
\begin{equation}
       \label{E5}
       p+q+r\ge
       -2\sqrt{qr} \, \cos \alpha \
       +2\sqrt{pr} \, \cos \beta  \
       +2\sqrt{pq} \, \cos \gamma.
\end{equation}
\end{statement}
\begin{proof}
Let us consider the quadratic function
\begin{equation}
       \label{E6}
       g\left(x\right)=
       x^2 - 2\left(\sqrt{r} \, \cos \beta
       +\sqrt{q} \, \cos \gamma \right)x
       + q + r + 2\sqrt{qr} \, \cos \alpha.
\end{equation}
Then a quarter of the discriminant is
\begin{equation}
        \label{E7}
        \frac{1}{4} \delta
        =\left(\sqrt{r} \, \cos \beta + \sqrt{q} \, \cos \gamma \right)^2
        -\left(q + r + 2\sqrt{qr} \, \cos \alpha \right).
\end{equation}
Based on $\alpha = \beta + \gamma $ we have $\cos \alpha = \cos  \left(\beta + \gamma \right)\
= \cos \beta \cos \gamma - \sin  \beta \sin \gamma$~and~hence
$$
\begin{array}{rcl}
\displaystyle\frac{1}{4} \delta
        &=& r \cos ^2 \beta
        + q \cos ^2 \gamma
        + 2 \sqrt{rq} \, \cos \beta \cos  \gamma
        - q - r - 2\sqrt{rq}\, \cos \alpha \\[1.5 ex]
        &=& r \cos ^2 \beta
        + q \cos ^2 \gamma
        + 2 \sqrt{rq} \, \cos \beta \cos \gamma
        - q - r - 2 \sqrt{rq} \, \cos \left(\beta + \gamma \right) \\[1.5 ex]
        &=&- r \sin ^2 \beta -q \sin ^2 \gamma
        + 2 \sqrt{rq} \, \cos \beta \cos \gamma
        - 2 \sqrt{rq} \, \cos \beta \cos \gamma
        + 2 \sqrt{rq} \, \sin \beta \sin \gamma.
\end{array}
$$
Using previous identity we obtained
$$
\delta
=
-4\left(\sqrt{r} \, \sin  \beta - \sqrt{q} \, \sin \gamma \right)^2<0,
$$
hence $g\left(x\right) \ge 0$.
Finally, letting $x=\sqrt{p}$ we obtained \eqref{E5}. \qed
\end{proof}
\begin{remark}
Let us emphasize that for term $A=p+q+r+2\sqrt{qr} \, \cos \alpha - 2\sqrt{pr}$
$\cos \beta -2\sqrt{pq} \, \cos \gamma$, when $\gamma= \alpha - \beta$,
follows inequality
$$
A
=
\left(\sqrt{r} -\sqrt{p} \, \cos  \beta + \sqrt{q} \, \cos \alpha \right)^2+\left(\sqrt{p} \, \sin  \beta - \sqrt{q} \, \sin \alpha \right)^2
\ge
0,
$$
analogously using the {\sc Lagrange}'s complete square identity from \cite{HojooLee}, \cite{Lu1}.
Therefore we have second proof of inequality (\ref{E5}).
\end{remark}
\begin{statement}
\label{St3}
Let $p, \, q, \, r \ge 0$ and $\alpha = \beta + \gamma$. Then we have the inequality
\begin{equation}
\label{E8}
       p+q+r\ge
       2\sqrt{qr}  \, \cos \alpha \
       -2\sqrt{pr} \, \cos \beta  \
       -2\sqrt{pq} \, \cos \gamma.
\end{equation}
\end{statement}
\begin{proof}
Let us consider the term $A = p+q+r-2\sqrt{qr} \, \cos \alpha +2\sqrt{pr} \, \cos \beta+2\sqrt{pq}$
$\cos \gamma$, for $\gamma = \alpha - \beta$. Notice that for the term $A$, by the {\sc Lagrange}'s
complete square identity, the following two representations are true.

\medskip
\noindent
$1^{\circ}\!\!\!.$ If $\displaystyle \frac{\pi}{2}\leq \alpha < \pi$, then $\cos \alpha \leq 0 $, and therefore

\vspace*{-4.0 mm}

\begin{equation}
\label{E9}
       A = \left(\sqrt{r} +\sqrt{p} \, \cos  \beta + \sqrt{q} \, \cos \alpha \right)^2
         + \left(\sqrt{p} \, \sin  \beta + \sqrt{q} \, \sin \alpha \right)^2
         - 4\sqrt{qr} \, \cos \alpha \ge 0.
\end{equation}

\noindent
$2^{\circ}\!\!\!\!.$ If $0 < \alpha <\displaystyle \frac{\pi}{2}$, then $\cos \alpha > 0$.
From $\alpha=\beta + \gamma $ follows $\cos \beta > 0$, and therefore

\vspace*{-4.0 mm}

\begin{equation}
\label{E10}
       A
       =
       \left(\sqrt{r} \!-\! \sqrt{p} \, \cos  \beta \!-\! \sqrt{q} \, \cos \alpha \right)^2
       \!+\!
       \left(\sqrt{p} \, \sin  \beta \!+\! \sqrt{q} \, \sin \alpha \right)^2
       \!+\!\,
       4\sqrt{pr} \, \cos \beta \ge 0.\!\!\!\!\!\!\!\!\qed
\end{equation}
\end{proof}

Let us introduce the division of the plane of triangle $\triangle\,ABC$ to following areas
$\lambda_0 = (+, +, +)$,
$\lambda_1 = (-, +, +)$,
$\lambda_2 = (+, -, +)$,
$\lambda_3 = (+, +, -)$,
$\lambda_4 = (+, -, -)$,
$\lambda_5 = (-, +, -)$,
$\lambda_6 = (-, -, +)$,
(Fig.$\;$2), via signes of homogenous barycentric coordinates of a point as given in the paper \cite{Yiu}
{\big (}see also the Section$\,$7.2$\,$in \cite{Farouki}{\big )}. Then~$\lambda_0$~is the interior area of
the triangle $\triangle\,ABC$. Let us notice that $(\lambda_0 \cup \lambda_1) \cup (BC)$ is the interior area of
the angle $\angle A$, and $\lambda_4$ is the interior area of the opposite angle. Analogously $(\lambda_0 \cup \lambda_2)
\cup (AC)$ is the interior area of the angle~$\angle B$, $\lambda_5$ is the interior area of the opposite angle
and $(\lambda_0 \cup \lambda_3) \cup (AB)$ is the interior area of the angle $\angle C$, $\lambda_6$
is the interior area of the opposite angle.

\newpage

\noindent
\begin{center} 

\vspace*{7.50 mm} \hspace*{-30.0 mm} \includegraphics*[height=30.0mm,keepaspectratio=true,scale=0.45]{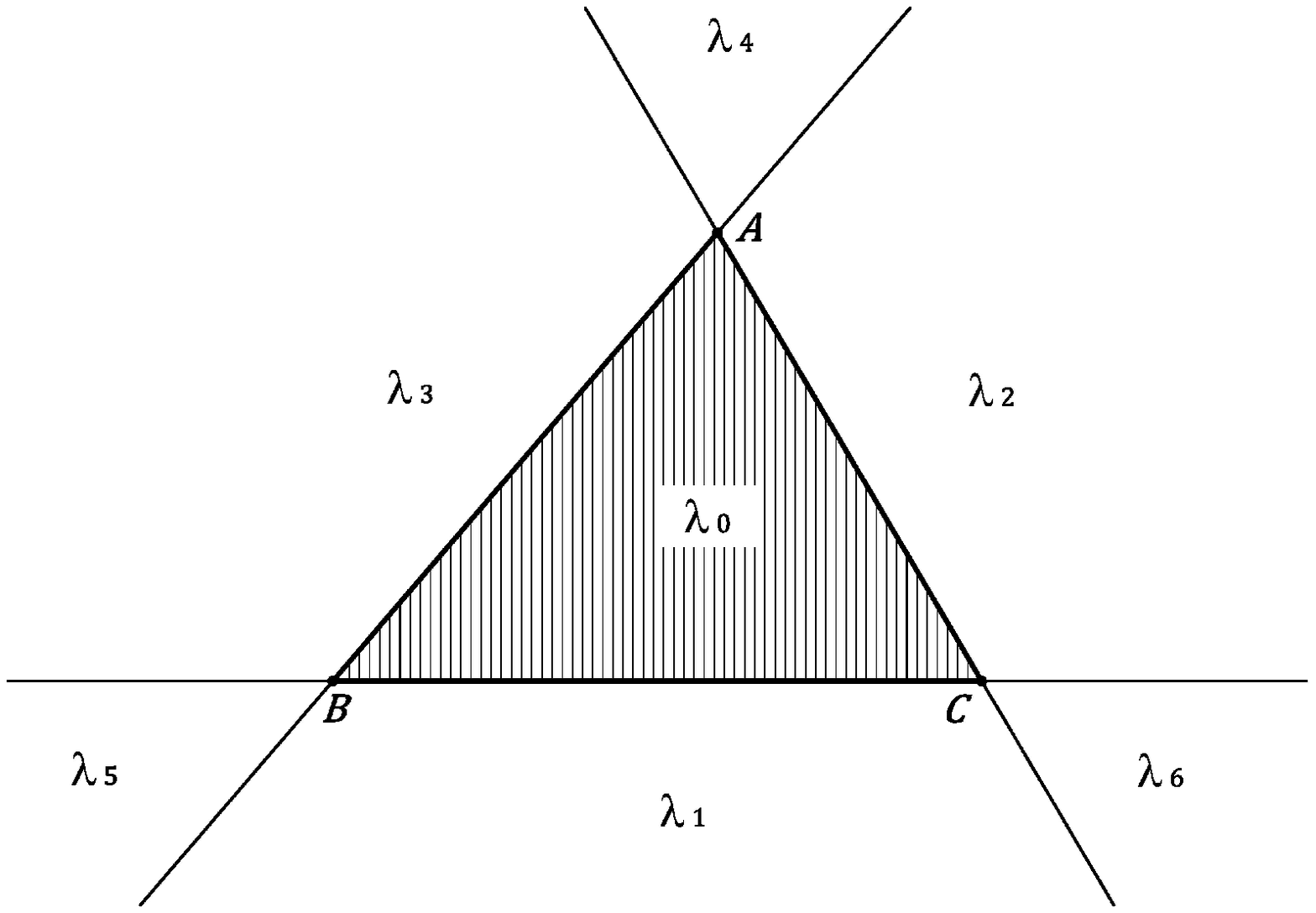}

\vspace*{0.00 mm}

\noindent
\textit{Figure 2: The division of the plane by}

\textit{the sidelines of the triangle $\triangle\,ABC$}

\end{center}

\vspace*{0.5 mm}

The following auxiliary statement is true.
\begin{lemma}
\label{Lemma0}
Let $B$ and $C$ be fixed points in the plane and let $M$ be arbitrary point in the plane.
For $\ell$ length of angle bisector of $\angle BMC$ from point $M$ following formulas are true$:$
\begin{equation}
\label{ell_1}
\ell
=
\displaystyle\frac{2 R_B R_C}{R_B + R_C} \, \cos \mbox{\small $\displaystyle\frac{\alpha_M}{2}$}
=
\displaystyle\frac{\sqrt{R_B R_C}}{R_B + R_C}\sqrt{(R_B + R_C)^2 - |BC|^2},
\end{equation}
where $R_B = |MB|$, $R_C = |MC|$ and $\alpha_M = \angle BMC$. Especially,
for $\pp{p}$ line throughout points $B$ and $C$ is true$:$

\vspace*{-4.5 mm}

\begin{equation}
\label{ell_0}
\ell
=
\left\{
\begin{array}{ccc}
0                                                             & : & M \in [BC],           \\[1.0 ex]
\mbox{\normalsize $\displaystyle\frac{2 R_B R_C}{R_B + R_C}$} & : & M \in \pp{p} \backslash [BC].
\end{array}
\right.
\end{equation}
\end{lemma}
\qquad In further considerations let $p = R_A$, $q = R_B$, $r = R_C$. Then, {\sc Z. Lu}, in the~paper \cite{Lu},
proved the following {\sc Barrow}'s type inequality.
\begin{theorem}
\label{Th0} \cite{Lu}
In the area $\lambda_0$ the following inequality is true$:$
\begin{equation}
\label{E11}
R_A \!+ R_B \!+ R_C \ge
    \!\left(\!\!\mbox{\small $\displaystyle\frac{\sqrt{R_C}}{\sqrt{R_B}}$}\!
    +\!\!\mbox{\small $\displaystyle\frac{\sqrt{R_B}}{\sqrt{R_C}}$}\right)\!\ell_a
    +\!\left(\!\!\mbox{\small $\displaystyle\frac{\sqrt{R_C}}{\sqrt{R_A}}$}\!
    +\!\!\mbox{\small $\displaystyle\frac{\sqrt{R_A}}{\sqrt{R_C}}$}\right)\!\ell_b
    +\!\left(\!\!\mbox{\small $\displaystyle\frac{\sqrt{R_A}}{\sqrt{R_B}}$}\!
    +\!\!\mbox{\small $\displaystyle\frac{\sqrt{R_B}}{\sqrt{R_A}}$}\right)\!\ell_c.
\end{equation}
\end{theorem}
\begin{remark}
{\sc Barrow}'s inequality is a consequence of the previous inequality.
\end{remark}
\qquad From previous Lemma follows next auxiliary statement.
\begin{lemma}
$\!\!$(i)$\;$ If $M=A$, i.e. $R_A=0$ then$:$
\begin{equation}
\label{E11a}
R_B + R_C
\geq
\left(
\!\!\mbox{\small $\displaystyle\frac{\sqrt{R_C}}{\sqrt{R_B}}$}
\!+\!
\mbox{\small $\displaystyle\frac{\sqrt{R_B}}{\sqrt{R_C}}$}\right)\!\ell_a.
\end{equation}

\vspace*{-2.5 mm}

\noindent
$\,$(ii)$\,$ If $M=B$, i.e. $R_B=0$ then$:$
\noindent
\begin{equation}
\label{E11b}
R_A + R_C
\geq
\left(
\!\!\mbox{\small $\displaystyle\frac{\sqrt{R_C}}{\sqrt{R_A}}$}
\!+\!
\mbox{\small $\displaystyle\frac{\sqrt{R_A}}{\sqrt{R_C}}$}\right)\!\ell_b.
\end{equation}

\vspace*{-2.5 mm}

\noindent
(iii) If $M=C$, i.e. $R_C=0$ then$:$
\begin{equation}
\label{E11c}
R_A + R_B
\geq
\left(
\!\!\mbox{\small $\displaystyle\frac{\sqrt{R_B}}{\sqrt{R_A}}$}
\!+\!
\mbox{\small $\displaystyle\frac{\sqrt{R_A}}{\sqrt{R_B}}$}\right)\!\ell_c.
\end{equation}
\end{lemma}

\break

Denote with \mbox{$\eu{cl}$} closure of a plane set. The following theorem~is~true.
\begin{theorem}
\label{Th1}
In the area $\mbox{$\eu{cl}$}\left(\lambda_1\right) \backslash \{B, C\}$ the following inequality is true$:$
\begin{equation}
\label{E12}
R_A \!+ R_B \!+ R_C \ge
    \!\!\left(\!\!\mbox{\small $\displaystyle\frac{\sqrt{R_C}}{\sqrt{R_B}}$}\!
    +\!\!\mbox{\small $\displaystyle\frac{\sqrt{R_B}}{\sqrt{R_C}}$}\right)\!\left(-\ell_a\right)
    +\!\left(\!\!\mbox{\small $\displaystyle\frac{\sqrt{R_C}}{\sqrt{R_A}}$}\!
    +\!\!\mbox{\small $\displaystyle\frac{\sqrt{R_A}}{\sqrt{R_C}}$}\right)\!\ell_b
    +\!\left(\!\!\mbox{\small $\displaystyle\frac{\sqrt{R_A}}{\sqrt{R_B}}$}\!
    +\!\!\mbox{\small $\displaystyle\frac{\sqrt{R_B}}{\sqrt{R_A}}$}\right)\!\ell_c.
\end{equation}
\end{theorem}
\begin{proof}
Let $M \!\in\! \mbox{$\eu{cl}$}\left(\lambda_1\right) \backslash \{B, C\}$, then $\alpha_M = \beta_M + \gamma_M$
i.e.$\;\displaystyle\frac{\alpha_M}{2} = \frac{\beta_M}{2} + \frac{\gamma_M}{2}$~(Fig.~3).

\begin{center} 

\vspace*{10.0mm} \hspace*{-18.0mm} \includegraphics*[height=35.0mm,keepaspectratio=true,scale=0.55]{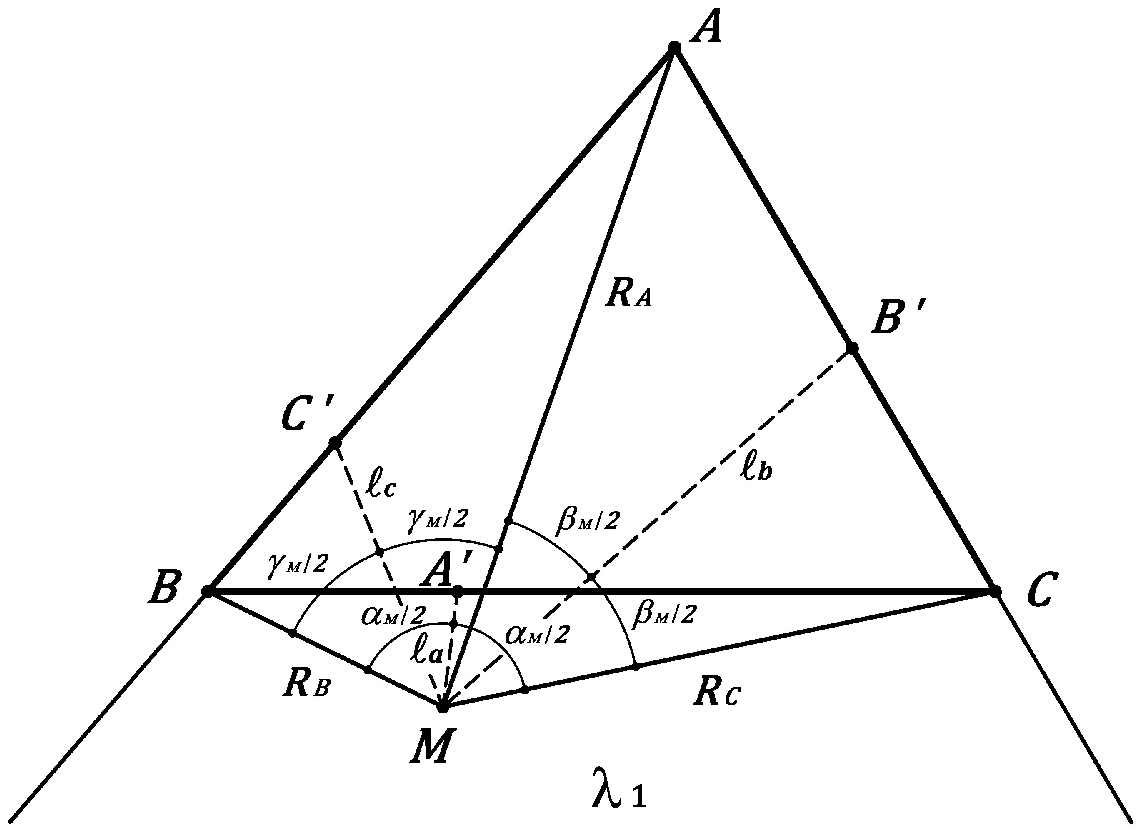}

\smallskip
\noindent
\textit{Figure 3:
Extension of the Barrow's inequality \\[1.0 ex]
for the point $M \in \mbox{$\eu{cl}$}\left(\lambda_1\right) \backslash \{B, C\}$}
\end{center}

\noindent
Based on the Statement~\ref{St2}, the following inequality holds
\begin{equation}
        \label{E14}
        R_A \!+\! R_B \!+\! R_C
        \ge
        - 2\sqrt{R_B R_C} \, \cos \frac{\alpha_M}{2}
        + 2\sqrt{R_A R_C} \, \cos \frac{\beta_M}{2}
        + 2\sqrt{R_A R_B} \, \cos \frac{\gamma_M}{2}
\end{equation}
and based on the Lemma \ref{Lemma0} from previous inequality we obtained
\begin{equation}
\label{E18}
\begin{array}{rcl}
\!\!\!\!\!
      R_A + R_B + R_C
      &\!\ge\!&
      \displaystyle - \frac{R_B + R_C}{\sqrt{R_B R_C}}\;\;\ell_a + \frac{R_A + R_C}{\sqrt{R_A R_C}}\;\; \ell_b + \frac{R_A + R_B}{\sqrt{R_A R_B}}\;\; \ell_c \\[3 ex]
      &\!=\!& \!\!\left(\!\!\mbox{\small $\displaystyle\frac{\sqrt{R_C}}{\sqrt{R_B}}$}\!
      + \!\mbox{\small $\displaystyle\frac{\sqrt{R_B}}{\sqrt{R_C}}$}\right)\!\left(-\ell_a\right)
      + \!\left(\!\!\mbox{\small $\displaystyle\frac{\sqrt{R_C}}{\sqrt{R_A}}$}\!
      + \!\mbox{\small $\displaystyle\frac{\sqrt{R_A}}{\sqrt{R_C}}$}\right)\!\ell_b
      + \!\left(\!\!\mbox{\small $\displaystyle\frac{\sqrt{R_A}}{\sqrt{R_B}}$}\!
      + \!\mbox{\small $\displaystyle\frac{\sqrt{R_B}}{\sqrt{R_A}}$}\right)\!\ell_c.\!\!\!\!\! \qed
      \end{array}
\end{equation}
\end{proof}

\medskip
Next two theorems are direct consequence of the Statement \ref{St2} by following cyclic replacements
$\alpha_M \mapsto \beta_M$, $\beta_M \mapsto \gamma_M$, $\gamma_M \mapsto \alpha_M$
and $R_A \mapsto R_B$, $R_B \mapsto R _C$, $R_C \mapsto R_A$ respectively.
\begin{theorem}
\label{Th2}
In the area $\mbox{$\eu{cl}$}\left(\lambda_2\right) \backslash \{A, C\}$ the following inequality is true$:$
\begin{equation}
\label{E19}
R_A \!+ R_B \!+ R_C \ge
    \!\!\left(\!\!\mbox{\small $\displaystyle\frac{\sqrt{R_C}}{\sqrt{R_B}}$}\!
    +\!\!\mbox{\small $\displaystyle\frac{\sqrt{R_B}}{\sqrt{R_C}}$}\right)\!\ell_a
    +\!\left(\!\!\mbox{\small $\displaystyle\frac{\sqrt{R_C}}{\sqrt{R_A}}$}\!
    +\!\!\mbox{\small $\displaystyle\frac{\sqrt{R_A}}{\sqrt{R_C}}$}\right)\!\!\left(-\ell_b\right)
    +\!\left(\!\!\mbox{\small $\displaystyle\frac{\sqrt{R_A}}{\sqrt{R_B}}$}\!
    +\!\!\mbox{\small $\displaystyle\frac{\sqrt{R_B}}{\sqrt{R_A}}$}\right)\!\!\ell_c.
\end{equation}
\end{theorem}
\begin{theorem}
\label{Th3}
In the area $\mbox{$\eu{cl}$}\left(\lambda_3\right) \backslash \{A, B\}$ the following inequality is true$:$
\begin{equation}
\label{E20}
R_A \!+ R_B \!+ R_C \ge
    \!\!\left(\!\!\mbox{\small $\displaystyle\frac{\sqrt{R_C}}{\sqrt{R_B}}$}\!
    \!+\!\!\mbox{\small $\displaystyle\frac{\sqrt{R_B}}{\sqrt{R_C}}$}\right)\!\ell_a
    +\!\left(\!\!\mbox{\small $\displaystyle\frac{\sqrt{R_C}}{\sqrt{R_A}}$}\!
    \!+\!\!\mbox{\small $\displaystyle\frac{\sqrt{R_A}}{\sqrt{R_C}}$}\right)\!\ell_b
    +\!\left(\!\!\mbox{\small $\displaystyle\frac{\sqrt{R_A}}{\sqrt{R_B}}$}\!
    \!+\!\!\mbox{\small $\displaystyle\frac{\sqrt{R_B}}{\sqrt{R_A}}$}\right)\!\!\left(-\ell_c\right).
\end{equation}
\end{theorem}

\medskip
The following theorem is true.

\begin{theorem}
\label{Th4}
In the area $\lambda_4$ the following inequality is true$:$
\begin{equation}
\label{E21}
R_A \!+\! R_B \!+\! R_C \ge
    \!\!\left(\!\!\mbox{\small $\displaystyle\frac{\sqrt{R_C}}{\sqrt{R_B}}$}\!
    \!+\!\!\mbox{\small $\displaystyle\frac{\sqrt{R_B}}{\sqrt{R_C}}$}\right)\!\ell_a
    \!+\!\left(\!\!\mbox{\small $\displaystyle\frac{\sqrt{R_C}}{\sqrt{R_A}}$}\!
    \!+\!\!\mbox{\small $\displaystyle\frac{\sqrt{R_A}}{\sqrt{R_C}}$}\right)\!\!\left(-\ell_b\right)
    \!+\!\left(\!\!\mbox{\small $\displaystyle\frac{\sqrt{R_A}}{\sqrt{R_B}}$}\!
    \!+\!\!\mbox{\small $\displaystyle\frac{\sqrt{R_B}}{\sqrt{R_A}}$}\right)\!\!\left(-\ell_c\right).
\end{equation}
\end{theorem}
\begin{proof}
Let $M \!\in\! \lambda_4$, then $\alpha_M=\beta_M+\gamma_M$ i.e.$\;\displaystyle\frac{\alpha_M}{2}
=\frac{\beta_M}{2}+\frac{\gamma_M}{2}$. Based on the Statement~\ref{St3} the following inequality is~true
\begin{equation}
\label{E22}
R_A + R_B + R_C \ge
        2\sqrt{R_B R_C} \, \cos \frac{\alpha_M}{2}
      - 2\sqrt{R_A R_C} \, \cos \frac{\beta_M}{2}
      - 2\sqrt{R_A R_B} \, \cos \frac{\gamma_M}{2}.
\end{equation}
Substitutions
\begin{equation}
\label{E15}
      \ell_a=\mid MA'\mid =2\frac{R_B R_C}{R_B + R_C} \cos \frac{\alpha_M}{2} \, ,
\end{equation}
\begin{equation}
\label{E16}
      \ell_b=\mid MB'\mid =2\frac{R_A R_C}{R_A + R_C} \cos \frac{\beta_M}{2} \, ,
\end{equation}
\begin{equation}
\label{E17}
      \ell_c=\mid MC'\mid =2\frac{R_A R_B}{R_A + R_B} \cos \frac{\gamma_M}{2} \,
\end{equation}
\noindent
in (\ref{E22}) give
\begin{equation}
\label{E23}
\begin{array}{rcl}
\!\!\!\!\!\!\!\!\!\!\!
        R_A \!+\! R_B \!+\! R_C
        &\!\!\ge\!&
        \displaystyle \frac{R_B + R_C}{\sqrt{R_B R_C}}\;\;\ell_a
        - \frac{R_A + R_C}{\sqrt{R_A R_C}}\;\; \ell_b
        - \frac{R_A + R_B}{\sqrt{R_A R_B}}\;\; \ell_c \\[3 ex]
        &\!\!=\!&\!\!\left(\!\!\mbox{\small $\displaystyle\frac{\sqrt{R_C}}{\sqrt{R_B}}$}\!
      \!+\!\,\!\mbox{\small $\displaystyle\frac{\sqrt{R_B}}{\sqrt{R_C}}$}\right)\!\ell_a
      \!+\!\,\!\left(\!\!\mbox{\small $\displaystyle\frac{\sqrt{R_C}}{\sqrt{R_A}}$}\!
      \!+\!\,\!\mbox{\small $\displaystyle\frac{\sqrt{R_A}}{\sqrt{R_C}}$}\right)\!\!\left(-\ell_b\right)
      \!+\!\,\!\left(\!\!\mbox{\small $\displaystyle\frac{\sqrt{R_A}}{\sqrt{R_B}}$}\!
      \!+\!\,\!\mbox{\small $\displaystyle\frac{\sqrt{R_B}}{\sqrt{R_A}}$}\right)\!\!\left(-\ell_c\right). \!\!\!\!\!\! \qed
    \end{array}
\end{equation}
\end{proof}

Next two theorems are direct consequence of the Statement \ref{St3} by following cyclic replacements
$\alpha_M \mapsto \beta_M$, $\beta_M \mapsto \gamma_M$, $\gamma_M \mapsto \alpha_M$
and $R_A \mapsto R_B$, $R_B \mapsto R _C$, $R_C \mapsto R_A$ respectively.

\begin{theorem}
\label{Th5}
In the area $\lambda_5$ the following inequality is true$:$
\begin{equation}
\label{E24}
R_A \!+\! R_B \!+\! R_C \ge
    \!\!\left(\!\!\mbox{\small $\displaystyle\frac{\sqrt{R_C}}{\sqrt{R_B}}$}\!
    \!+\!\!\mbox{\small $\displaystyle\frac{\sqrt{R_B}}{\sqrt{R_C}}$}\right)\!\!\left(-\ell_a\right)
    \!+\!\left(\!\!\mbox{\small $\displaystyle\frac{\sqrt{R_C}}{\sqrt{R_A}}$}\!
    \!+\!\!\mbox{\small $\displaystyle\frac{\sqrt{R_A}}{\sqrt{R_C}}$}\right)\!\ell_b
    \!+\!\left(\!\!\mbox{\small $\displaystyle\frac{\sqrt{R_A}}{\sqrt{R_B}}$}\!
    \!+\!\!\mbox{\small $\displaystyle\frac{\sqrt{R_B}}{\sqrt{R_A}}$}\right)\!\!\left(-\ell_c\right).
\end{equation}
\end{theorem}

\begin{theorem}
\label{Th6}
In the area $\lambda_6$ the following inequality is true$:$
\begin{equation}
\label{E25}
R_A \!+\! R_B \!+\! R_C \ge
    \!\!\left(\!\!\mbox{\small $\displaystyle\frac{\sqrt{R_C}}{\sqrt{R_B}}$}\!
    \!+\!\!\mbox{\small $\displaystyle\frac{\sqrt{R_B}}{\sqrt{R_C}}$}\right)\!\!\left(-\ell_a\right)
    \!+\!\left(\!\!\mbox{\small $\displaystyle\frac{\sqrt{R_C}}{\sqrt{R_A}}$}\!
    \!+\!\!\mbox{\small $\displaystyle\frac{\sqrt{R_A}}{\sqrt{R_C}}$}\right)\!\!\left(-\ell_b\right)
    \!+\!\left(\!\!\mbox{\small $\displaystyle\frac{\sqrt{R_A}}{\sqrt{R_B}}$}\!
    \!+\!\!\mbox{\small $\displaystyle\frac{\sqrt{R_B}}{\sqrt{R_A}}$}\right)\!\ell_c.
\end{equation}
\end{theorem}

Now, we give definition of the signed angle bisector for the point $M$ in the plane of the triangle $\triangle\,ABC$.
Let be $A$ fixed vertex and let $\pp{p}$ be line through vertices $B$~and~$C$. Denote $d\!=\!|MA_1|$ distance of the
point $M$ to the line $\pp{p}$ and let $\ell\!=\!|MA'|$ be length of the bisector of the angle $\angle BMC$.
If $d'$ be signed distance of the point $M$ to the line $\pp{p}$ related to the vertex $A$ \cite{Saul} (p.$\;$308.),
then $d'\!=\!+d$ if $M$ and $A$ with same side of line $\pp{p}$, otherwise $d'\!=\!-d$. Let us define signed angle
bisector $\ell'$ analogously $\ell'\!=\!+\ell$ if $M$ and $A$ with same side of line $\pp{p}$,
otherwise $\ell'\!=\!-\ell$ (Fig. 4). In the case $M \!\in\! \pp{p}$ then~$d'\!=\!0$ and
then $\ell'$ given by formula (\ref{ell_0}).

\begin{center} 

\vspace*{+13.0 mm} \hspace*{-25.0 mm}
\includegraphics*[height=35.0mm,keepaspectratio=true,scale=0.55]{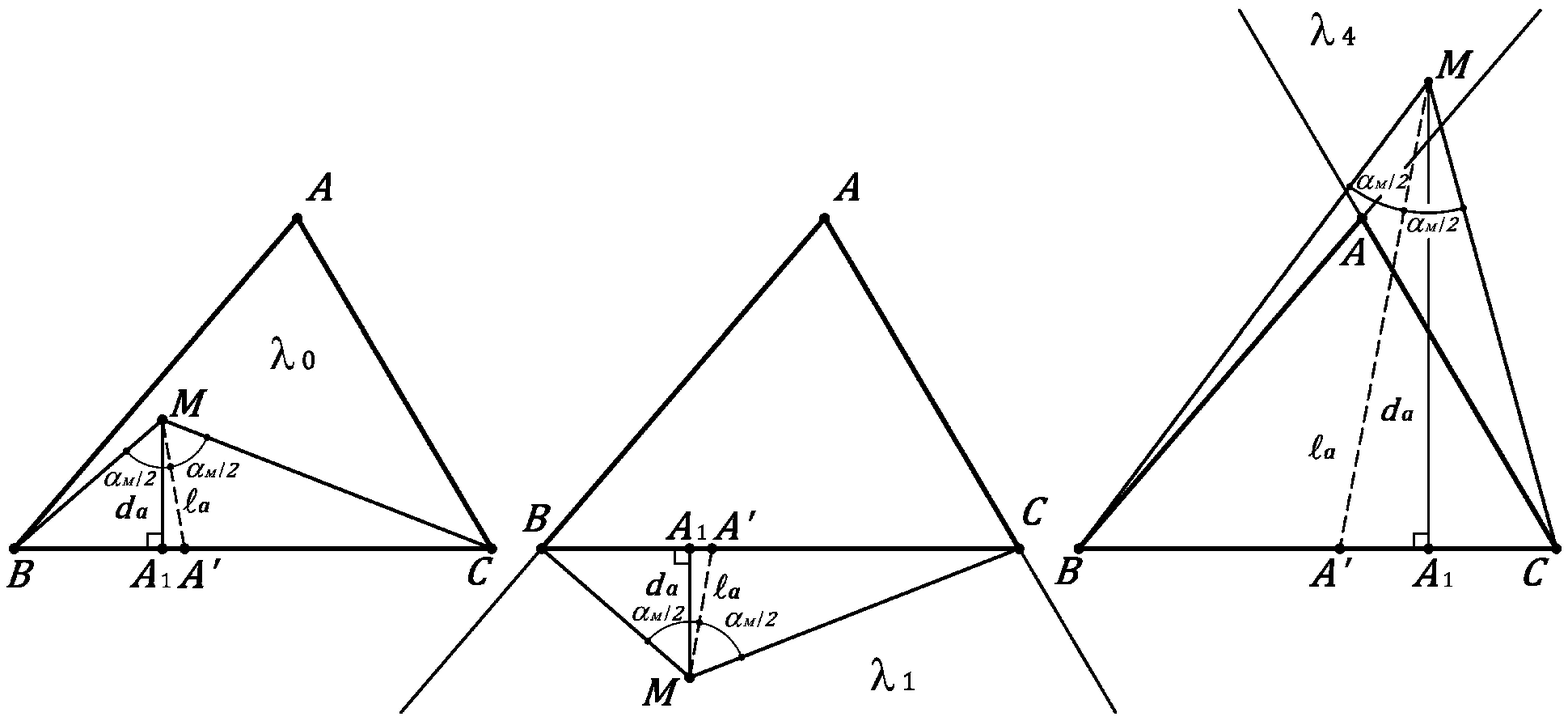}
{\small
$$
\,\,\,\,\,\,\,\,\,\,\,
d'_a=\!+d_a  \; \; \; \; \; \; \; \; \; \; \; \; \; \; \; \; \; \; \; \; \; \; \; \; \; \; \; \; \; \; \; \; \; \; \; \; \; \; \; \;
d'_a=-d_a    \; \; \; \; \; \; \; \; \; \; \; \; \; \; \; \; \; \; \; \; \; \; \; \; \; \; \; \; \; \; \; \; \; \; \; \; \; \; \; \; \;
d'_a=+d_a
$$
\vspace*{-4.5 mm}
$$
\,\,\,\,\,\,\,\,\,\;
\ell'_a=+\,\ell_a   \; \; \; \; \; \; \; \; \; \; \; \; \; \; \; \; \; \; \; \; \; \; \; \; \; \; \; \; \; \; \; \; \; \; \; \; \; \; \; \;
\ell'_a=-\,\ell_a   \; \; \; \; \; \; \; \; \; \; \; \; \; \; \; \; \; \; \; \; \; \; \; \; \; \; \; \; \; \; \; \; \; \; \; \; \; \; \; \; \; \;
\ell'_a=\!+\,\ell_a
$$}

\vspace*{-2.5 mm}

\noindent
\textit{Figure 4: Signed distances and signed angle bisectors}

\end{center}

\medskip
Let us denote
$
\mu_{1} \!\!=\!\! \mbox{$\eu{cl}$}(\lambda_1) \backslash \{B,C\}\!,\;
\mu_{2} \!\!=\!\! \mbox{$\eu{cl}$}(\lambda_2) \backslash \{A,C\}\!,\;
\mu_{3} \!\!=\!\! \mbox{$\eu{cl}$}(\lambda_3) \backslash \{A,B\}\!,\;
\mu_{4} \!\!=\!\! \lambda_{4},\;
\mu_{5} \!\!=\!\! \lambda_{5}
$
and
$
\mu_{6} \!\!=\!\! \lambda_{6}
$.
Then $\bigcup_{\,i=1}^{\,6}{\mu_{i}} \cup \{A, B, C\}$ is a complete division of the plane
of the triangle $\triangle\,ABC\!$. Finally, analogously to {\sc Dergiades} extension
of the {\sc Erd\" os}- {\sc Mordell} inequal\-ity \cite{Derg}, from previous theorems,
an extension of {\sc Barrow}'s type in\-equality (\ref{E11}) is obtained by the following theorem.

\vspace*{-1.5 mm}

\begin{statement}
\label{St4}
For the point $M \!\in\! \displaystyle\bigcup_{i=1}^{6}{\,\mu_{i}}$ the following inequality is true$:$

\vspace*{-2.5 mm}

\begin{equation}
\label{E26}
\;\;\;\;\;\;\;\;\;
 R_A + R_B + R_C
    \ge
        \left(\!\!\mbox{\small $\displaystyle\frac{\sqrt{R_C}}{\sqrt{R_B}}$}\!
    \!+\!\mbox{\small $\displaystyle\frac{\sqrt{R_B}}{\sqrt{R_C}}$}\right)\!\!\ell_a'
    \!+\!\left(\!\!\mbox{\small $\displaystyle\frac{\sqrt{R_C}}{\sqrt{R_A}}$}\!
    \!+\!\mbox{\small $\displaystyle\frac{\sqrt{R_A}}{\sqrt{R_C}}$}\right)\!\!\ell_b'
    \!+\!\left(\!\!\mbox{\small $\displaystyle\frac{\sqrt{R_A}}{\sqrt{R_B}}$}\!
    \!+\!\mbox{\small $\displaystyle\frac{\sqrt{R_B}}{\sqrt{R_A}}$}\right)\!\!\ell_c'\,;
\end{equation}
otherwise for points $M \!=\! A$, $M \!=\! B$, $M \!=\! C$ following inequalities
$(\ref{E11a})$, $(\ref{E11b})$, $(\ref{E11c})$ are true respectively.
\end{statement}


\end{document}